\documentclass[12pt,a4paper]{amsart}
\usepackage{amsmath,amssymb,amsthm,bm,bbm}
\usepackage{array,mathdots,blkarray,arydshln}
\usepackage{graphicx,enumerate}
\usepackage{layout}

\theoremstyle{plain}
\newtheorem{theorem}{Theorem}[section]
\newtheorem{lemma}[theorem]{Lemma}

\theoremstyle{definition}
\newtheorem{remark}[theorem]{Remark}

\newtheorem{question}[theorem]{Question}
\newtheorem{example}[theorem]{Example}
\newtheorem*{acknowledgments}{Acknowledgments}

\numberwithin{equation}{section}

\newcommand{\bQ}{\mathbbm{Q}}

\newcommand{\bP}{\mathbbm{P}}
\newcommand{\bA}{\mathbbm{A}}

\title[A two-dimensional rationality problem and intersections of two quadrics]
{A two-dimensional rationality problem and intersections of two quadrics}

\author[A. Hoshi]{Akinari Hoshi}
\address{Department of Mathematics,
Niigata University, Niigata 950-2181, Japan}
\email{hoshi@math.sc.niigata-u.ac.jp}

\author[M. Kang]{Ming-chang Kang}
\address{Department of Mathematics,
National Taiwan University, Taipei, Taiwan}
\email{kang@math.ntu.edu.tw}

\author[H. Kitayama]{Hidetaka Kitayama}
\address{Department of Mathematics,
Wakayama University, Wakayama 640-8510, Japan}
\email{hkitayam@wakayama-u.ac.jp}

\author[A. Yamasaki]{Aiichi Yamasaki}
\address{Department of Mathematics,
Kyoto University, Kyoto 606-8502, Japan}
\email{aiichi.yamasaki@gmail.com}

\thanks{
This work was partially supported by JSPS KAKENHI Grant Number 25400027,
15K17511, 16K05059, 19K03418.
Parts of the work were finished when the
first-named author, the third-named author and the fourth-named author
were visiting the National Center for Theoretic Sciences (Taipei),
whose support is gratefully acknowledged.}


\begin{document}
\begin{abstract}
Let $k$ be a field with char $k\neq 2$ and $k$ be not algebraically closed. 
Let $a\in k\setminus k^2$ and $L=k(\sqrt{a})(x,y)$ be a field extension of 
$k$ where $x,y$ are algebraically independent over $k$. 
Assume that $\sigma$ is a $k$-automorphism on $L$ defined by 
\[
\sigma: \sqrt{a}\mapsto -\sqrt{a},\ 
x\mapsto \frac{b}{x},\ y\mapsto \frac{c(x+\frac{b}{x})+d}{y}
\]
where $b,c,d \in k$, $b\neq 0$ and at least one of $c,d$ is non-zero.
Let $L^{\langle\sigma\rangle}=\{u\in L:\sigma(u)=u\}$ be the fixed subfield 
of $L$. 
We show that $L^{\langle\sigma\rangle}$ is isomorphic to the function field 
of a certain surface in $\bP^4_k$ which is given as the intersection of 
two quadrics. 
We give criteria for the $k$-rationality of $L^{\langle\sigma\rangle}$ 
by using the Hilbert symbol. 
As an appendix of the paper, we also give an alternative geometric proof 
of a part of the result which is provided to the authors by 
J.-L. Colliot-Th\'el\`ene. 
\end{abstract}

\maketitle

\section{Introduction} \label{se1}

Let $k$ be a field and $L$ be a finitely generated extension field
of $k$. $L$ is called {\it $k$-rational} (or {\it rational over
$k$}) if $L$ is purely transcendental over $k$, i.e.\ $L$ is
$k$-isomorphic to the quotient field of some polynomial ring over $k$.
$L$ is called {\it stably $k$-rational} if $L(y_1,\ldots,y_m)$ is
$k$-rational for some $y_1,\ldots,y_m$ which are algebraically
independent over $L$.
$L$ is called {\it $k$-unirational} if $L$ is $k$-isomorphic to a
subfield of some $k$-rational field.
It is obvious that ``$k$-rational" $\Rightarrow$ ``stably $k$-rational"
$\Rightarrow$ ``$k$-unirational".
The L\"uroth problem asks, under what situation, the converse is true,
i.e. ``$k$-unirational" $\Rightarrow$ ``$k$-rational".
The L\"uroth problem is a famous problem in algebraic geometry.
For a survey of it and related rationality problems
(e.g. Noether's problem), see, for examples, \cite{MT} and \cite{Sw}.

Throughout this paper, we assume, unless otherwise specified, 
that $k$ is a field with char $k\neq 2$ 
and $k$ is not an algebraically closed field.

\begin{question}\label{Q1}
Let $a\in k\setminus k^2$ and $L=k(\sqrt{a})(x,y)$ be a
field extension of $k$ where $x,y$ are algebraically independent
over $k(\sqrt{a})$.
Let $\sigma$ be the $k$-automorphism on $L$ defined by
\[
\sigma: \sqrt{a}\mapsto -\sqrt{a},\
x\mapsto x,\ y\mapsto \frac{f(x)}{y}
\]
where $f(x) \in k[x]$ is a non-zero polynomial.
Let $L^{\langle\sigma\rangle}=\{u\in L:\sigma(u)=u\}$
be the fixed subfield of $L$.
When will $L^{\langle\sigma\rangle}$ be $k$-rational?
\end{question}

Hajja, Kang and Ohm \cite{HKO} gave an answer to Question \ref{Q1} 
when the degree of $f(x)$ is $\le 2$: 

\begin{theorem}[{\cite[Theorem 6.7]{HKO}, \cite[Theorem 4.2]{Ka2}, see also \cite[Lemma 4.1]{HKK} for the last statement}]
\label{t1.2}
Let $k$ be a field with {\rm char} $k\neq 2$
and $k(\sqrt{a})$ be a separable quadratic field extension of $k$
with $a \in k$.
Write ${\rm Gal}(k(\sqrt{a})/k)= \langle \sigma \rangle$ and
extend the action of $\sigma$ to $k(\sqrt{a})(x,y)$ by
\[
\sigma: \sqrt{a} \mapsto -\sqrt{a},\ x \mapsto x,\ y\mapsto
\frac{f(x)}{y}\quad (f(x)\in k[x]\setminus \{0\}).
\]
{\rm (1)} When $f(x)=b$, $k(\sqrt{a})(x,y)^{\langle\sigma\rangle}$ is
$k$-rational if and only if the Hilbert symbol $(a,b)_{2,k}=0$.\\
{\rm (2)} When $\deg f(x)=1$,
$k(\sqrt{a})(x,y)^{\langle\sigma\rangle}$ is always $k$-rational.\\
{\rm (3)} When $f(x)=b(x^2-c)$ for some $b,c\in k\setminus \{0\}$, then
$k(\sqrt{a})(x,y)^{\langle\sigma\rangle}$ is $k$-rational if and only
if $(a,b)_{2,k}\in {\rm Br}(k(\sqrt{ac})/k)$, i.e.
$(a,b)_{2,k(\sqrt{ac})}=0$. 

Moreover, if $k(\sqrt{a})(x,y)^{\langle\sigma\rangle}$ is not
$k$-rational, then $k$ is an infinite field,
${\rm Br}(k)$ is non-trivial, and
$k(\sqrt{a})(x,y)^{\langle\sigma\rangle}$ is not $k$-unirational.
\end{theorem}

In the general case, 
Question \ref{Q1} is related to the rationality of conic bundles 
of $\bP^1$ over a non-closed field $k$ 
investigated by Iskovskikh \cite{Is1, Is2, Is3}; 
also see Lemma \ref{l2.1} and Yamasaki \cite{Ya2}. 


The purpose of this paper is to find a solution of 
the following similar question:

\begin{question}\label{Q2}
Let $a\in k\setminus k^2$ and $L=k(\sqrt{a})(x,y)$ be a
field extension of $k$ where $x,y$ are algebraically independent
over $k(\sqrt{a})$.
Let $\sigma$ be the $k$-automorphism on $L$ defined by
\[
\sigma: \sqrt{a}\mapsto -\sqrt{a},\ x\mapsto \frac{b}{x},\
y\mapsto \frac{c(x+\frac{b}{x})+d}{y}
\]
where $b,c,d \in k$, $b\neq 0$ and at least one of $c,d$ is non-zero.

Let $L^{\langle\sigma\rangle}=\{u\in L:\sigma(u)=u\}$
be the fixed subfield of $L$. 
Find a necessary and sufficient condition 
for the field $L^{\langle\sigma\rangle}$ to be $k$-rational.
\end{question}

\begin{remark}
We may consider the following variant of Question \ref{Q2}. 
Let $L=k(x,y)$ be a field extension of $k$ 
where $x,y$ are algebraically independent over $k$.
Let $\sigma$ be the $k$-automorphism on $L$ defined by
\[
\sigma: x\mapsto \frac{b}{x},\ y\mapsto \frac{c(x+\frac{b}{x})+d}{y}
\]
where $b,c,d \in k$, $b\neq 0$ and at least one of $c,d$ is 
non-zero.

Find a necessary and sufficient condition 
for the field $L^{\langle\sigma\rangle}$ to be $k$-rational. 
Be aware that $\sigma$ acts trivially on the ground field $k$ of
the field $L$, while in Question \ref{Q2} $\sigma$ acts faithfully 
on $k(\sqrt{a})$ which is the ground field of $L$. 
In this case, Kang \cite[Theorem 2.4]{Ka1} proved that 
$L^{\langle\sigma\rangle}$ is always $k$-rational. 
This type of the action appears in the related rationality problems 
naturally (see e.g. \cite{HKY}, \cite{Ya1}, \cite{HK}). 
\end{remark}

Now return to Question \ref{Q2}. 
To begin with, we describe explicitly 
the fixed field $L^{\langle\sigma\rangle}$. 
It turns out that 
the fixed field $L^{\langle\sigma\rangle}$ 
is isomorphic to the function field of a 
certain surface in $\bP^4_k$ which is 
given as the intersection of two quadrics. 

\begin{lemma}\label{l1.5}
Let the notation be the same as in Question \ref{Q2}.
Then the fixed field $L^{\langle\sigma\rangle}$ is $k$-isomorphic
to the field $k(t_i:1\leq i\leq 4)$ with the relations
$t_1^2-at_2^2=b$ and $t_3^2-at_4^2=2ct_1+d$.

Conversely, if $L_0=k(t_i:1\leq i\leq 4)$ is a field extension of
$k$ satisfying that ${\rm trdeg}_k(L_0)=2$, $t_1^2-at_2^2=b$ and
$t_3^2-at_4^2=2ct_1+d$, then $L_0$ is $k$-isomorphic to the fixed
subfield $k(\sqrt{a})(x,y)^{\langle\sigma\rangle}$ where the action
of $\sigma$ is given in Question \ref{Q2}.
\end{lemma}

\begin{proof}
Given the action of $\sigma$ on $k(\sqrt{a})(x,y)$, define
\begin{align*}
t_1&=\frac{1}{2}\left(x+\frac{b}{x}\right),\
t_2=\frac{1}{2\sqrt{a}}\left(x-\frac{b}{x}\right),\\
t_3&=\frac{1}{2}\left(y+\frac{c(x+\frac{b}{x})+d}{y}\right),\
t_4=\frac{1}{2\sqrt{a}}\left(y-\frac{c(x+\frac{b}{x})+d}{y}\right).
\end{align*}
It is routine to verify that
$t_i \in L^{\langle\sigma\rangle}$ for $1 \le i \le 4$.
In fact, we will show that
$L^{\langle\sigma\rangle}=k(t_i:1\leq i\leq 4)$
and these generators satisfy the relations
\begin{align*}
\begin{cases}
t_1^2-at_2^2=b,\\
t_3^2-at_4^2=2ct_1+d.
\end{cases}
\end{align*}

First note that $k(t_i:1\leq i\leq 4)\subset L^{\langle\sigma\rangle}
\subset L$.

We claim that $k(\sqrt{a})(t_i:1\leq i\leq 4)=L$.
The elements $x+\frac{b}{x}$ and $x-\frac{b}{x}$ belong to the field
$k(\sqrt{a})(t_i:1\leq i\leq 4)$.
Thus $x$ belongs to $k(\sqrt{a})(t_i:1\leq i\leq 4)$.
Similarly, $y$ belongs to $k(\sqrt{a})(t_i:1\leq i\leq 4)$.
It follows that $k(\sqrt{a})(t_i:1\leq i\leq 4)=L$.

Since $[L:k(t_i:1\leq i\leq 4)]
=[k(\sqrt{a})(t_i:1\leq i\leq 4):k(t_i:1\leq i\leq 4)]\leq 2$
and $[L:L^{\langle\sigma\rangle}]=2$,
we find that $L^{\langle\sigma\rangle}=k(t_i:1\leq i\leq 4)$.

On the other hand, if $L_0$ is the field extension of $k$ given
in the statement of this lemma,
write $x= t_1+\sqrt{a}t_2$, $y=t_3+\sqrt{a}t_4$.
It is easy to verify the assertion.
\end{proof}

Because of Lemma \ref{l1.5}, 
we are led to the following question 
which is equivalent to Question \ref{Q2}. 

\begin{question}\label{Q4}
Let $k$ be a field with {\rm char} $k\neq 2$,
and $a,b,c,d\in k$ satisfying that $a\in k\setminus k^2$, $b\neq 0$
and at least one of $c$ and $d$ is non-zero.
Define the projective surface $X$ in $\bP^4_{k}$ 
by the intersection of the following two quadrics: 
\begin{align*}
\begin{cases}
f_1:T_1^2-a T_2^2-bT_5^2=0,\\
f_2:T_3^2-a T_4 ^2-2c T_1T_5-dT_5^2=0
\end{cases}
\end{align*}
and take the affine open subset $X_0\subset X$ defined by $T_5 \neq 0$. 
Explicitly, $X_0$ is the closed subset of the affine $4$-space 
$\bA^4_k$ defined by 
\begin{align}
\begin{cases}
t_1^2-at_2^2=b,\\
t_3^2-at_4^2=2ct_1+d
\end{cases}\label{eq11}
\end{align}
where $t_i:=T_i/T_5$ for $1 \le i \le 4$. 
We take the function field $k(X_0)$ of $X_0$ over $k$. 

Find a necessary and sufficient condition for 
the function field $k(X_0)$ to be $k$-rational. 

Clearly $k(X_0)$ (or $k(X)$) is isomorphic to $L_0$ over $k$ 
(for the definition of $L_0$, see Lemma \ref{l1.5}). 
\end{question}
\begin{remark}\label{r1.6}
An answer to Question \ref{Q4} can be determined up to sign of $c$
because $L_0$ is stable under the action of $t_1\mapsto -t_1$ and
$2c(-t_1)+d=2(-c)t_1+d$.
\end{remark}

The main result of this paper is the following theorem 
which gives an answer to Question \ref{Q2} and Question \ref{Q4}. 
We will write $(a,b)_{2,k}$ for the Hilbert symbol over $k$,
i.e. the quaternion algebra over $k$,
where $a, b$ are non-zero elements of $k$
(see \cite[Section 11]{Dr}).
By the definition, $(a,b)_{2,k}=0$ if and only if there are elements
$x,y\in k$ such that $ax^2+by^2=1$.
\begin{theorem}\label{t1.8}
Let $X$ be the intersection of two quadrics 
$($the projective surface$)$ in $\bP^4_{k}$ as in Question \ref{Q4}. 
Then the variety $X$ has all of its $k$-points in the affine open subset 
$X_0$ defined by $T_5\neq 0$. 
The open subset $X_0$ has at most one singular $k$-point, 
which exists if and only if $d^2-4bc^2=0$, and it corresponds to 
the $k$-point $(-\tfrac{d}{2c},0,0,0)$. 
And the following statements are equivalent:\\
{\rm (i)} The field $L_0=k(X_0)$ is $k$-rational;\\
{\rm (ii)} The field $L_0=k(X_0)$ is stably $k$-rational;\\
{\rm (iii)} The field $L_0=k(X_0)$ is $k$-unirational;\\
{\rm (iv)} The surface $X_0$ $($or $X$$)$ has at least two $k$-points;\\
{\rm (v)} The surface $X_0$ $($or $X$$)$ has a smooth $k$-point. 

When $d^2-4bc^2\neq 0$ $($i.e. when $X_0$ is smooth$)$, 
the statements {\rm (i)}--{\rm (v)} are also equivalent to 
each of the following:\\
{\rm (vi)} There exist elements $\alpha, \beta \in k$ such that 
$\alpha ^2-a\beta^2=b$, and either $2c\alpha+d=0$, or 
$2c\alpha+d\neq 0$ with $(a,2c\alpha+d)_{2,k(\sqrt{d^2-4bc^2})}=0$;\\  
{\rm (vii)} $(a,b)_{2,k}=0$, and 
for any elements $\alpha,\beta\in k$ such that $\alpha^2-a\beta^2=b$, 
we have either $2c\alpha + d =0$, or
$2c\alpha+d\neq 0$ with $(a,2c\alpha+d)_{2,k(\sqrt{d^2-4bc^2})}=0$.

When $d^2-4bc^2=0$ $($i.e. when $X_0$ is singular$)$, 
the statements {\rm (i)}--{\rm (v)} are also equivalent to\\
{\rm (viii)} $(a,2d)_{2,k}=0$. 

Moreover, if $L_0$ is not
$k$-rational, then $k$ is an infinite field and
${\rm Br}(k)$ is non-trivial.

Furthermore, the conditions are also equivalent 
after changing $c$ to $-c$.
\end{theorem}

\begin{remark}
The last statement of Theorem \ref{t1.8} 
for {\rm (iv)}, {\rm (v)}, {\rm (vi)} and {\rm (vii)} 
yields non-trivial phenomenon although it is trivial 
for {\rm (i)}, {\rm (ii)}, {\rm (iii)} and {\rm (viii)} 
(see Remark \ref{r1.6}). 
\end{remark}

We organize this paper as follows.
In Section \ref{se2}, the proof of Theorem \ref{t1.8} will be given. 
Section \ref{se3} contains some applications and examples 
of Theorem \ref{t1.8}. 
As an appendix of the paper, in Section \ref{se4}, 
we give an alternative geometric proof of a part of Theorem \ref{t1.8} 
which is provided to the authors by J.-L. Colliot-Th\'el\`ene. 

\begin{acknowledgments}
We thank J.-L. Colliot-Th\'el\`ene who called our attention to a 
rationality criterion in the paper of Colliot-Th\'el\`ene, 
Sansuc and Swinnerton-Dyer \cite{CTSSD}. 
He was so nice to explain to one of us (A. H.) how to deduce 
the $k$-rationality in the presence of a smooth $k$-point 
(see Section $\ref{se4}$). 

We are grateful to the referee who supplied a very detailed report. 
This paper is greatly improved by his/her helpful suggestions. 
In particular, he/she pointed out an erroneous argument of 
the proof of Theorem \ref{t1.8} 
in a previous version of this paper. 
\end{acknowledgments}

\section{Proof of Theorem \ref{t1.8}}\label{se2}

By Theorem \ref{t1.2}, we may obtain the following lemmas. 
The proof of the following lemma is omitted because
it is almost the same as that of Lemma \ref{l1.5}.

\begin{lemma}\label{l2.1}
Let the notation be the same as in Theorem \ref{t1.2}.
Then the fixed field $k(\sqrt{a})(x,y)^{\langle\sigma\rangle}$ is
$k$-isomorphic to the field $k(x,t_1,t_2)$ with the relation
$t_1^2-at_2^2=f(x)$.

Conversely, if $L_0=k(x,t_1,t_2)$ is a field extension of $k$
satisfying that ${\rm trdeg}_k(L_0)=2$ and  $t_1^2-at_2^2=f(x)$,
a non-zero polynomial in $k[x]$, then $L_0$ is $k$-isomorphic to
the fixed subfield $k(\sqrt{a})(x,y)^{\langle\sigma\rangle}$ where
the action of $\sigma$ is given in Theorem \ref{t1.2}.
\end{lemma}

\begin{lemma}\label{l2.2}
Let $k$ be a field with {\rm char} $k\neq 2$, and $a, b, c \in k$ with
$ab \neq 0$. Let $L_0=k(x, t_1, t_2)$ be a field extension of $k$
such that ${\rm trdeg}_k(L_0)=2$ and $t_1 ^2 - a t_2 ^2 = b(x^2 -c)$.
Then the following statements are equivalent:\\
{\rm (i)} $L_0$ is $k$-rational;\\
{\rm (ii)} $L_0$ is $k$-unirational;\\
{\rm (iii)} $(a,b)_{2,k(\sqrt{ac})}=0$.
\end{lemma}

\begin{proof}
By Lemma \ref{l2.1}, $L_0$ is $k$-isomorphic to the fixed subfield
$k(\sqrt{a})(x,y)^{\langle \sigma \rangle}$ of the field
$k(\sqrt{a})(x,y)$ where $\sigma: \sqrt{a} \mapsto -\sqrt{a},\ x
\mapsto x,\ y\mapsto \frac{b(x^2-c)}{y}$.
Thus we may apply Theorem \ref{t1.2}.

If $c\neq 0$, apply Theorem \ref{t1.2} {\rm (3)} and we are done.

Suppose that $c=0$. Then $k(\sqrt{ac})=k$ and 
we will show that $L_0$ is $k$-rational 
if and only if $(a,b)_{2,k}=0$. 
Since $L_0=k(x,t_1,t_2)=k(\frac{t_1}{x},\frac{t_2}{x},x)$
with relation $(\frac{t_1}{x})^2-a(\frac{t_2}{x})^2=b$,
we may apply Theorem \ref{t1.2} {\rm (1)}. 
The equivalence ${\rm (i)}\Leftrightarrow {\rm (ii)}$ 
follows from the last statement of Theorem \ref{t1.2}. 
\end{proof}

The verification of the following lemma is routine. We may use the
multiplicative property of the norm from the field $k(\sqrt{a})$ to
the field $k$ to get an alternative (heuristic) proof.

\begin{lemma}\label{l2.4}
Let $a, x_1, x_2, y_1, y_2$ be elements in some field.
Then $(x_1^2 -a x_2^2)(y_1^2 -a y_2^2)=(x_1y_1+ax_2y_2)^2-a(x_1y_2+x_2y_1)^2$.
If $y_1^2 -a y_2^2 \neq 0$, we have also
$\frac{x_1^2 -a x_2^2}{y_1^2 -a y_2^2}
=\left(\frac{x_1y_1-ax_2y_2}{y_1^2-ay_2^2}\right)^2
-a\left(\frac{x_1y_2-x_2y_1}{y_1^2-ay_2^2}\right)^2$.
\end{lemma}

\bigskip
{\it Proof of Theorem \ref{t1.8}.}

First we find that 
for a $k$-point $[t_1:t_2:t_3:t_4:t_5]$ on $X$, 
it is necessary that $t_5\neq 0$.
Otherwise, we have $t_1^2-a t_2^2=0$ and $t_3^2-a t_4 ^2=0$ 
where $t_i \in k$ for $1\le i \le 4$. 
If some $t_i$ is non-zero, we find that $a \in k^2$. 
This is impossible from the assumption. 
In short, all $k$-points on $X$ belong to the affine open subset 
$X_0$ defined by $T_5 \neq 0$. 

We may find the singular locus of $X_0$ 
by applying the Jacobian criterion for the affine varieties 
\cite[page 130, Theorem 2.19]{Liu}. It is not difficult to find all the
singularities on $X_0$ by considering the Jacobian matrix
\begin{align*}
\left[\frac{\partial f_i}{\partial t_j}\right]_{1\leq i\leq 2,
1\leq j\leq 4}=
\left[
\begin{array}{ccccc}
2 t_1&-2at_2&0&0\\
-2c&0&2t_3&-2at_4
\end{array}\right].
\end{align*}

We find that either 
(i) $X_0$ has no singular points, 
or 
(ii) the point $(-\tfrac{d}{2c},0,0,0)$ 
is the unique singular point with a necessary 
condition that $d^2-4bc^2=0$. 

Namely, the surface $X_0$ is smooth if and only if $d^2-4bc^2\neq 0$, 
and if $d^2-4bc^2=0$, then there exists 
the unique singular point $(-\tfrac{d}{2c},0,0,0)$ on $X_0$. 
The first statement of Theorem \ref{t1.8} follows. 

\bigskip
{\it Proof of ${\rm (iv)}\Leftrightarrow {\rm (v)}$ of Theorem \ref{t1.8}.} ---------

{}From the argument above, we see that 
${\rm (iv)}\Rightarrow {\rm (v)}$. 
The reverse implication 
${\rm (v)}\Rightarrow {\rm (iv)}$ 
also follows whenever $d^2-4bc^2=0$. 
If $d^2-4bc^2\neq 0$ and there exists 
a smooth $k$-point $(\alpha,\beta,\gamma,\delta)$ on $X_0$ 
which satisfies 
$\alpha^2-a\beta^2=b$, 
$\gamma^2-a\delta^2=2c\alpha+d$, 
then we find that $(\beta,\gamma,\delta)\neq (0,0,0)$. 
Suppose $(\beta,\gamma,\delta)=(0,0,0)$. 
Then we have $\alpha^2=b$ and $2c\alpha+d=0$,  
and hence $d^2-4bc^2=(d+2\alpha c)(d-2\alpha c)=0$, 
which is a contradiction. 
Thus there exist at least two $k$-points 
$(\alpha,\pm\beta,\pm\gamma,\pm\delta)$ on $X_0$. 
Namely, we obtain ${\rm (v)}\Rightarrow {\rm (iv)}$ 
when $d^2-4bc^2\neq 0$ also. \qed\\

The implications ${\rm (i)}\Rightarrow {\rm (ii)}\Rightarrow {\rm (iii)}$ 
of Theorem \ref{t1.8} are standard. 
It remains to show that 
${\rm (iii)}\Rightarrow {\rm (v)}\Rightarrow {\rm (i)}$, 
${\rm (i)}\Leftrightarrow {\rm (vi)}$ and 
${\rm (i)}\Leftrightarrow {\rm (vii)}$ 
under the assumption $d^2-4bc^2\neq 0$, 
and ${\rm (i)}\Leftrightarrow {\rm (viii)}$ 
under the assumption $d^2-4bc^2=0$. 
After that, the last statement follows from Remark \ref{r1.6}. 

\bigskip
{\it Proof of ${\rm (iii)}\Rightarrow {\rm (v)}$ of Theorem \ref{t1.8}.} ---------

Assume that $k$ is a finite field. 
Then 
we have $(a,b)_{2,k}=0$ for any $a,b\in k\setminus\{0\}$. 
Hence we may take $\alpha,\beta\in k$ with $\alpha^2-a\beta^2=b$ 
and also $\gamma,\delta\in k$ with $\gamma^2-a\delta^2=2c\alpha+d$ 
whenever $2c\alpha+d\neq 0$. 
Because we have $(\gamma,\delta)\neq (0,0)$, 
under the assumption $2c\alpha+d\neq 0$, 
we obtain a smooth $k$-point 
$(\alpha,\beta,\gamma,\delta)$ on $X_0$ from the argument above. 
Suppose that $2c\alpha+d=0$. 
If $d^2-4bc^2\neq 0$, then we just 
take a $k$-point $(\alpha,\beta,0,0)$ on $X_0$ 
which is a smooth $k$-point from the argument above. 
If $d^2-4bc^2=0$, then $(\beta,\gamma,\delta)=(0,0,0)$ may occur.
But in this case, we have $\alpha\neq 0$, $cd\neq 0$
and we can take $\alpha^\prime:=-\alpha$ such that
$\alpha^{\prime\,2}-a\beta^2=b$,
$2c\alpha^\prime+d\neq 0$ (remember that char $k\neq 2$).
This implies that the existence of a smooth $k$-point 
$(\alpha^\prime,\beta,\gamma^\prime,\delta^\prime)$ on $X_0$ 
such that $\alpha^{\prime\,2}-a\beta^2=b$,
$\gamma^{\prime\,2}-a\delta^{\prime\,2}=2c\alpha^\prime+d$
with $(\beta,\gamma^\prime,\delta^\prime)\neq (0,0,0)$.

Assume that $k$ is an infinite field.\\

Step 1.
Usually a $k$-unirational field $L$ is defined as a field $L$ with
an embedding $L\rightarrow k(x_1,\ldots,x_n)$ where
$k(x_1,\ldots,x_n)$ is a rational function field with $n$ variables over $k$
and $n$ is some positive integer.
However, a stronger result is possible by Ohm \cite{Oh}: 
It is shown that,
if $k$ is any field and $L$ is $k$-unirational,
then there is a $k$-embedding $L\rightarrow k(x_1,\ldots,x_n)$ where
$k(x_1,\ldots,x_n)$ is a rational function field and $n$ variables over $k$
where $n={\rm trdeg}_k(L)$.

Thus, from the assumption that $L_0$ is $k$-unirational,
we may assume that $L_0$ is a subfield of $k(x,y)$ where $x$ and $y$
are algebraically independent over $k$.
It follows that every element of $L_0$ can be written as a rational function
of $x$ and $y$.
In particular, there exist polynomials $h(x,y)$, $h_i(x,y)\in k[x,y]$
such that $h(x,y)$ is not zero and
\begin{align}
t_i=\frac{h_i(x,y)}{h(x,y)}\label{eq21}
\end{align}
for $1\leq i\leq 4$. 
Not all of $h_2(x,y), h_3(x,y), h_4(x,y)$ are the zero polynomial;
otherwise, $L_0=k(t_1,t_2,t_3,t_4)$ becomes an extension of $k$ 
with transcendence degree at most one. 
We may assume that ${\rm gcd}\{h(x,y),h_i(x,y):1\leq i\leq 4\}=1$.
Substitute (\ref{eq21}) into (\ref{eq11}).
We obtain
\begin{align*}
\begin{cases}
\displaystyle{
\frac{h_1(x,y)^2}{h(x,y)^2}-a\frac{h_2(x,y)^2}{h(x,y)^2}}&=b,\\
\displaystyle{
\frac{h_3(x,y)^2}{h(x,y)^2}-a\frac{h_4(x,y)^2}{h(x,y)^2}}&=
\displaystyle{
2c\frac{h_1(x,y)}{h(x,y)}+d}.
\end{cases}
\end{align*}

\bigskip
Step 2.
Since $k$ is an infinite field, 
we may find infinitely many $a_0, b_0 \in k$ with 
$h(a_0,b_0)h_j(a_0,b_0) \neq 0$ 
where $h(x,y)h_j(x,y)$ is not the zero polynomial for some $2\leq j\leq 4$. 
Then we have 
\begin{align*}
\begin{cases}
\alpha^2-a\beta^2=b,\\
\gamma^2-a\delta^2=2c\alpha+d
\end{cases}
\end{align*}
with $(\beta,\gamma,\delta)\neq (0,0,0)$. 
We get a smooth $k$-point $(\alpha,\beta,\gamma,\delta)$ on $X_0$ 
from the argument above and hence (v) follows.\qed

\bigskip
{\it Proof of ${\rm (v)}\Rightarrow {\rm (i)}$ of Theorem \ref{t1.8}.} ---------

We take a smooth $k$-point $(\alpha,\beta,\gamma,\delta)$ on $X_0$ 
which satisfies 
$\alpha^2-a\beta^2=b$, 
$\gamma^2-a\delta^2=2c\alpha+d$ 
and $(\alpha,\beta,\gamma,\delta)\neq (-\tfrac{d}{2c},0,0,0)$. 
Then, by (\ref{eq11}), we obtain
$(t_1^2-\alpha^2)-a(t_2^2-\beta^2)=0$.
Rewrite this relation as
\begin{align*}
\frac{t_1-\alpha}{t_2+\beta}=a\frac{t_2-\beta}{t_1+\alpha}.
\end{align*}

Define $t_0=\frac{t_2-\beta}{t_1+\alpha}$.

We claim that $k(t_1,t_2)=k(t_0)$.
Besides the condition that
$t_0=\frac{t_2-\beta}{t_1+\alpha}$,
we have also $at_0=\frac{t_1-\alpha}{t_2+\beta}$.
These two conditions are equivalent to
\begin{align*}
\begin{cases}
t_0t_1-t_2&\!\!\!\!\!=-\alpha t_0-\beta,\\
t_1-at_0t_2&\!\!\!\!\!=\alpha+a\beta t_0.
\end{cases}
\end{align*}

Solve the above linear simultaneous equation of $t_1$ and $t_2$.
We find
\begin{align*}
\begin{cases}
(1-at_0^2)t_1=\alpha(1+at_0^2)+2a\beta t_0,\\
(1-at_0^2)t_2=\beta(1+at_0^2)+2\alpha t_0.
\end{cases}
\end{align*}
Thus $k(t_1,t_2)=k(t_0)$ as we expect. It follows that we have that
$L_0=k(t_0,t_3,t_4)$ with the relation
\begin{align}
(t_3^2-at_4^2)(1-at_0^2)=a(2c\alpha-d)t_0^2+4ac\beta t_0+(2c\alpha+d).
\label{eq23}
\end{align}
By Lemma \ref{l2.4}, we obtain the identity that
\[
(t_3^2-at_4^2)(1-at_0^2)=(t_3+at_0t_4)^2-a(t_0t_3+t_4)^2.
\]
Thus we define
\[
t_5=t_3+at_0t_4,\ t_6=t_0t_3+t_4.
\]
Then we have $L_0=k(t_0,t_5,t_6)$ and the relation provided by
Equation ({\ref{eq23}) becomes
\begin{align*}
t_5^2-at_6^2=a(2c\alpha-d)t_0^2+4ac\beta t_0+(2c\alpha+d).
\end{align*}

Define $T_0=\frac{1}{t_0}$,
$T_5=\frac{t_5}{t_0}$, $T_6=\frac{t_6}{t_0}$.
Then we have $L_0=k(t_0,t_5,t_6)=k(T_0,T_5,T_6)$ and
\begin{align}
T_5^2-aT_6^2=(2c\alpha+d)T_0^2+4ac\beta T_0+a(2c\alpha-d).\label{eq24}
\end{align}
We will solve the rationality problem according to the cases 
$2c\alpha+d\neq 0$ and $2c\alpha+d=0$.\\

Case 1. Assume that $2c\alpha+d\neq 0$.

Define $t_7=T_0+\frac{2ac\beta}{2c\alpha+d}$.
Then by Equation (\ref{eq24})
we have $L_0=k(T_5,T_6,t_7)$ and
\begin{align}
T_5^2-aT_6^2=(2c\alpha+d)t_7^2-\frac{a(d^2-4bc^2)}{2c\alpha+d}.\label{eq25}
\end{align}

It follows from the assumption that 
$\gamma^2-a\delta^2=2c\alpha+d\neq 0$. 
Define $t_8=\frac{T_5\gamma-aT_6\delta}{\gamma^2-a\delta^2},
t_9=\frac{T_5\delta-T_6\gamma}{\gamma^2-a\delta^2}$. 
Then $L_0=k(t_7,t_8,t_9)$ and it follows from Lemma \ref{l2.4} that 
\begin{align*}
t_8^2-at_9^2=t_7^2-\frac{a(d^2-4bc^2)}{(2c\alpha+d)^2}. 
\end{align*} 

Define $T_7=(2c\alpha+d)t_7$, $T_8=(2c\alpha+d)t_8$, 
$T_9=(2c\alpha+d)t_9$.
Then we get that $L_0=k(T_7,T_8,T_9)$ and
\begin{align}
T_8^2-aT_9^2=T_7^2-a(d^2-4bc^2).\label{eq26}
\end{align}

We may apply Lemma \ref{l2.2} to Equation (\ref{eq26}) (with $b=1$
where $b$ is a element in the notation of Lemma \ref{l2.2}).
Since the Hilbert symbol $(a,1)_{2,k}=0$ (as $b=1$ is a square in $k$),
we find that
$(a,1)_{2,k(\sqrt{d^2-4bc^2})}=0$. By Lemma \ref{l2.2},
$L_0=k(T_7, T_8, T_9)$ is $k$-rational.\\

Case 2. Assume that $2c\alpha+d=0$. 

Then we have $c\neq 0$ because we assume that 
at least one of $c$ and $d$ is non-zero. 

First we show that $d^2-4bc^2\neq 0$. 

Suppose $d^2-4bc^2=0$. 
By the assumption $2c\alpha+d=0$, we have 
$d^2-4c^2\alpha^2=(d+2c\alpha)(d-2c\alpha)=0$ 
and hence $b=\alpha^2$ (note that $c\neq 0$). 
Substitute $b=\alpha^2$ into $\alpha^2-a \beta^2=b$, 
obtain $\beta=0$. 
On the other hand, we have 
$\gamma^2-a\delta^2=2c\alpha+d=0$. 
Since $\sqrt{a} \notin k$, it follows that $\gamma=\delta=0$. 
Hence we obtain that $(\beta,\gamma,\delta)=(0,0,0)$. 
This contradicts to the assumption 
$(\alpha,\beta,\gamma,\delta)\neq (-\tfrac{d}{2c},0,0,0)$. 

Next we show that $\beta \neq 0$. 

Suppose $\beta=0$. 
Then we find that $b=\alpha^2$.
Hence $d^2-4bc^2=(d+2c\alpha)(d-2c\alpha)=0$.
This is also impossible because we have $d^2-4bc^2\neq 0$.

By the assumption $2c\alpha+d = 0$, 
Equation (\ref{eq24}) becomes 
$T_5^2-aT_6^2=4ac\beta T_0+a(2c\alpha - d)$.

Thus we find that $T_0 \in k(T_5,T_6)$ because 
$4ac\beta \neq 0$ and hence 
$L_0=k(T_0, T_5, T_6)=k(T_5,T_6)$ is $k$-rational.\qed

\bigskip
{\it Proof of ${\rm (i)}\Rightarrow {\rm (vi)}$ of Theorem \ref{t1.8}.} ---------

We have the assumption $d^2-4bc^2\neq 0$. 

Since ${\rm (i)}\Rightarrow {\rm (v)}$, 
there exists 
a smooth $k$-point $(\alpha,\beta,\gamma,\delta)$ on $X_0$ 
which satisfies 
$\alpha^2-a\beta^2=b$, 
$\gamma^2-a\delta^2=2c\alpha+d$. 
If $2c\alpha+d=0$, then nothing should be done. 
If $2c\alpha+d \neq 0$, then 
we have $(a,b)_{2,k}=0$ and $(a,2c\alpha+d)_{2,k}=0$.
Hence $(a,2c\alpha+d)_{2,k(\sqrt{d^2-4bc^2})}=0$ follows.\qed

\bigskip
{\it Proof of ${\rm (vi)}\Rightarrow {\rm (i)}$ of Theorem \ref{t1.8}.} ---------

We have the assumption $d^2-4bc^2\neq 0$. 

The condition $\alpha^2-a\beta^2=b$ guarantees that the part of
${\rm (v)}\Rightarrow {\rm (i)}$ (till Equation (\ref{eq24})) is still valid
in the present situation. 
Thus $L_0=k(T_0,T_5,T_6)$ with a relation defined
by Equation (\ref{eq24}).

We consider the cases $2c\alpha+d \neq 0$ and $2c\alpha+d=0$ separately.\\

Case 1. Assume that $2c\alpha+d\neq 0$.

We have $L_0=k(T_5,T_6,t_7)$ with the relation defined by Equation (\ref{eq25}).
By the assumption that the Hilbert symbol
$(a,2c\alpha+d)_{2,k(\sqrt{d^2-4bc^2})}=0$, 
we find that $L_0$ is $k$-rational by Lemma \ref{l2.2}.\\

Case 2. Assume that $2c\alpha+d=0$.

Equation (\ref{eq24}) becomes $T_5^2-aT_6^2=4ac\beta T_0+a(2c\alpha-d)$. 
We find that $c\beta\neq 0$. 
This follows from Case 2 of the proof of 
${\rm (v)}\Rightarrow {\rm (i)}$ of Theorem \ref{t1.8} 
(now we have the assumption $d^2-4bc^2\neq 0$).

Because $4ac\beta\neq 0$, we have $T_0 \in k(T_5, T_6)$ and 
$L_0=k(T_5,T_6)$ is $k$-rational.
\qed

\bigskip
{\it Proof of ${\rm (i)}\Leftrightarrow {\rm (vii)}$ of Theorem \ref{t1.8}.} ---------

We have the assumption $d^2-4bc^2\neq 0$. 

By ${\rm (vi)}$, it is obvious that $(a, b)_{2,k}=0$ is a necessary
condition of ${\rm (i)}$. It remains to show that, if $\alpha, \beta \in k$
satisfying that $\alpha ^2-a \beta ^2=b$, then $L_0$ is $k$-rational
if and only if either $2c\alpha + d =0$, or $2c\alpha + d \neq 0$
with $(a, 2c\alpha + d)_{2, k(\sqrt{d^2-4bc^2})}=0$. 

By Lemma \ref{l1.5}, we will assume that $L_0= L^{\langle\sigma\rangle}$
where $L=k(\sqrt{a})(x,y)$ and $\sigma$ acts on it by the $k$-automorphism
\[
\sigma: \sqrt{a}\mapsto -\sqrt{a},\ x\mapsto \frac{b}{x},\ y\mapsto
\frac{c(x+\frac{b}{x})+d}{y}
\]
where $b=\alpha^2-a\beta^2$ with $\alpha,\beta\in k$. 

Define $z=\sqrt{a}\left(\frac{1+s}{1-s}\right)$ 
where $s=\frac{x}{\alpha-\sqrt{a}\beta}$ and 
$w=y(z+\sqrt{a})$. 
It is not difficult to verify that
$k(\sqrt{a})(x,y)= k(\sqrt{a})(z,w)$ and
\begin{align*}
\sigma: \sqrt{a}\mapsto -\sqrt{a},\ z \mapsto z,\ w \mapsto \frac{f(z)}{w}
\end{align*}
where $f(z)=(2c\alpha+d)z^2+4ac\beta z +a(2c\alpha-d)$.\\

Case 1. Assume that $2c\alpha+d\neq 0$.

We will apply Lemma \ref{l2.1} and Lemma \ref{l2.2} because we may write 
\[
f(z)= (2c\alpha+d)\left[\left(z+ \frac{2ac\beta}{2c\alpha+d}\right)^2
-\frac{a(d^2- 4bc^2)}{(2c\alpha+ d)^2}\right].
\]

By Lemma \ref{l2.1} and Lemma \ref{l2.2}, the fixed field 
$k(\sqrt{a})(z,w)^{\langle\sigma\rangle}$ is $k$-rational 
if and only if $(a,2c\alpha+d)_{2,k^{\prime}}=0$ where 
$k^{\prime}=k(\sqrt{a\frac{a(d^2-4bc^2)}{(2c\alpha+d)^2}})
=k(\sqrt{d^2-4bc^2})$. 
Done.\\

Case 2. Assume that $2c\alpha+d=0$.

It is necessary that $c \neq 0$ because $(c,d)\neq (0,0)$.
Suppose that $\beta=0$.  
Then $b=\alpha^2$ and $d^2-4bc^2=(d+2c\alpha)(d-2c\alpha)=0$. 
But this is impossible by the assumption $d^2-4bc^2\neq 0$. 
Hence we have $c\beta\neq 0$. 

In short, $f(z)=4ac\beta z + a(2c\alpha-d)$ is a degree one polynomial 
in $z$. 
We may apply Theorem \ref{t1.2} ${\rm (2)}$ to conclude that 
$k(\sqrt{a})(z,w)^{\langle\sigma\rangle}$ is $k$-rational 
without any extra condition.\qed
%

\bigskip
{\it Proof of ${\rm (i)}\Leftrightarrow {\rm (viii)}$ of Theorem \ref{t1.8}.} ---------

We have the assumption $d^2-4bc^2=0$. 
Then $c\neq 0$ and $b=(\frac{d}{2c})^2\in k^2$. 

By Lemma \ref{l1.5}, we may assume that $L_0= L^{\langle\sigma\rangle}$
where $L=k(\sqrt{a})(x,y)$ and $\sigma$ acts on it by the $k$-automorphism
\[
\sigma: \sqrt{a}\mapsto -\sqrt{a},\ x\mapsto \frac{b}{x},\ y\mapsto
\frac{c(x+\frac{b}{x})+d}{y}
\]
where $b=(\frac{d}{2c})^2\in k^2$.
Define $\beta^\prime=\frac{d}{2c}$ and $u, v$ by
\begin{align*}
u=\, \sqrt{a}\frac{\beta^\prime-x}{\beta^\prime+x},\,
v= \,(\sqrt{a}- u)y.
\end{align*}
One may verify that $L=k(\sqrt{a})(x,y)=k(\sqrt{a})(u,v)$ and $\sigma$
acts on $u,v$ by
\[
\sigma: \sqrt{a}\mapsto -\sqrt{a},\ u\mapsto u,\
v\mapsto \frac{-2ad}{v}.
\]
It follows from Lemma \ref{l2.1} that 
$L_0=k(u,V_1,V_2)$ with $V_1^2-aV_2^2=-2ad$. 

By Theorem \ref{t1.2} ${\rm (1)}$,
$L_0$ is $k$-rational if and only if the symbol $(a, -2ad)_{2,k}=0$.
Hence the assertion follows
because $(a,-2ad)_{2,k}=(a,-a)_{2,k}+(a,2d)_{2,k}=(a,2d)_{2,k}$.\qed

\section{Applications and examples}\label{se3}

We thank J.-L. Colliot-Th\'el\`ene who pointed out that 
the algebraic manipulations in the proof of 
Theorem \ref{t1.8} in Section \ref{se2} yield 
the following theorem. 
\begin{theorem}\label{t3.1}
Let $X$ be the intersection of two quadrics 
$($the projective surface$)$ in $\bP^4_k$ 
and $X_0\subset X$ be the affine open subset defined by $T_5 \neq 0$ 
as in Question \ref{Q4}.\\
{\rm (1)} Assume that $d^2-4bc^2\neq 0$ $($i.e. when $X_0$ is smooth$)$.\\
{\rm (1-1)} If $X$ has a $k$-point $[\alpha:\beta:\gamma:\delta:1]$ 
with $2c\alpha+d\neq 0$, then 
$X$ is birationally equivalent over $k$ 
to the surface in $\bP^3$ 
defined by the quadric $T_1^2-aT_2^2-T_3^2+a(d^2-4bc^2)T_4^2=0$;\\
{\rm (1-2)} If $X$ has a $k$-point $[\alpha:\beta:\gamma:\delta:1]$ 
with $2c\alpha+d=0$, then 
$X$ is birationally equivalent over $k$ 
to the surface in $\bP^3$ 
defined by the quadric $T_5^2-aT_6^2-T_7T_8=0$;\\
{\rm (2)} Assume that $d^2-4bc^2=0$ $($i.e. when $X_0$ is singular$)$. 
Then $X$ is stably birationally equivalent over $k$ 
to the curve in $\bP^2$ defined by the quadric $V_1^2-aV_2^2+2adV_3^2=0$. 
\end{theorem}
\begin{proof}
(1-1) follows from Equation (\ref{eq26}) 
in Case $1$ of the proof of 
${\rm (v)}\Rightarrow {\rm (i)}$ of Theorem \ref{t1.8}.\\
(1-2) follows from Case $2$ in the proof of 
${\rm (v)}\Rightarrow {\rm (i)}$ of Theorem \ref{t1.8} by 
defining $T_7=4ac\beta T_0+a(2c\alpha-d)$.\\ 
(2) follows from the proof of 
${\rm (i)}\Leftrightarrow {\rm (viii)}$ of Theorem \ref{t1.8}.
\end{proof}

We also get the following theorem as an application of Theorem \ref{t1.8}.
\begin{theorem}\label{t3.2}
Let $k$ be a field with {\rm char} $k\neq 2$, 
and $X$ be an irreducible projective variety over $k$ 
of dimension $2n$ $($$n$ is a positive integer$)$. 
Assume that $X$ is birationally equivalent over $k$ to the product 
of the projective surfaces $X_i$ $($$1 \le i \le n$$)$ 
which are defined in the analogous way as in Question \ref{Q4} 
and the parameters $a_i, b_i, c_i, d_i$ depend on the surface $X_i$. 
If $X$ is $k$-unirational, then $X$ is $k$-rational.
\end{theorem}
\begin{proof}
Suppose that $X$ is $k$-unirational. Then each $X_i$ is also $k$-unirational. 
Apply Theorem \ref{t1.8}. 
We find that each $X_i$ is $k$-rational. 
Hence $X$ is $k$-rational.
\end{proof}

We give examples of Theorem \ref{t1.8} 
for the case where $X_0$ is smooth  
and also the case where $X_0$ is singular. 

\begin{example}[Examples of Theorem \ref{t1.8}]\label{e3.3}
Let the notation be the same as in Theorem \ref{t1.8}. 
Then $X$ is the intersection of two quadrics 
(the projective surface) in $\bP^4_k$ and 
$X_0\subset X$ is the affine open subset defined by $T_5 \neq 0$ 
as in Question \ref{Q4} 
and $L_0$ is the function field of $X$ (or $X_0$) over $k$ 
as in Lemma \ref{l1.5}. 
We will give examples of Theorem \ref{t1.8} with $k=\bQ$.\\

(1) 
The case where $X_0$ is smooth, i.e. $d^2-4bc^2\neq 0$. 

We take $k=\bQ$ and 
$(a,b,c,d,\alpha,\beta)= (2,1,c,c,-1,0)$ with 
$\alpha^2-a\beta^2=b$ and $c\neq 0$. 

We will find a necessary and sufficient condition such that 
$L_0$ is $\bQ$-rational.

Note that $2c\alpha+d=-c \neq 0$ and $d^2-4bc^2=-3c^2\neq 0$. 
By Theorem \ref{t1.8} (vi), 
we find that $L_0$ is $\bQ$-rational 
$\Leftrightarrow$ $(2,-c)_{2,\bQ(\sqrt{-3})}=0$.

(i) If $a$ and $1-a$ are non-zero elements of $k$, from
$a\cdot 1^2+(1-a)\cdot 1^2=1$, we find that $(a,1-a)_{2,k}=0$. 
Put $a=2$ into $(a,1-a)_{2,k}=0$. 
We find that $(2, -1)_{2,k}=0$ for any field $k$. 
Thus $(2,-c)_{2,k(\sqrt{-3})}=(2,c)_{2,k(\sqrt{-3})}$. 

(ii) For the special case that $c=p$ is a prime number,
we know that $(2,2)_{2,\bQ}=(2,-1)_{2,\bQ}=0$ and 
for $p\geq 3$, 
by the Hasse principle, 
$(2,p)_{2,\bQ}=0$ 
$\Leftrightarrow$ $(2,p)_{2,\bQ_l}=0$ for any prime $l\leq\infty$ 
$\Leftrightarrow$ $(2,p)_{2,\bQ_l}=0$ for $l=2,p$ 
$\Leftrightarrow$ the Legendre symbol $(\frac{2}{p})=1$ 
$\Leftrightarrow$ $p\equiv\pm 1\pmod{8}$. 
Consequently, if $c$ is a prime number $p$ and 
$(a,b,c,d)= (2,1,p,p)$, then $L_0$ is $\bQ$-rational 
provided that $p=2$ or $p\equiv\pm 1\pmod{8}$ $(p\geq 3)$. 

On the other hand, over $K=\bQ(\sqrt{-3})$,
we see that $(2,3)_{2,\bQ}\neq 0$ 
but $(2,3)_{2,K}=0$ because $3$ is ramified in $K$. 
By the Hasse principle over $K$, 
if $c$ is a prime number $p\geq 5$ with $p\equiv\pm 3\pmod{8}$, 
then $(2,p)_{2,K}= 0\Leftrightarrow $
$(2,p)_{2,K_\mathfrak{p}}=0$ for any prime $\mathfrak{p}$ in $K$ 
$\Leftrightarrow$ $(2,p)_{2,K_\mathfrak{p}}=0$ for $\mathfrak{p}$ 
lying above $p$ (and $2$) 
$\Leftrightarrow$ $p$ is inert in $K$ 
$\Leftrightarrow$ the Legendre symbol $(\frac{-3}{p})=-1$ 
$\Leftrightarrow$ $p\equiv 2\pmod{3}$.

In summary, if $c=p$ is a prime number and $(a,b,c,d)= (2,1,p,p)$,
$L_0$ is $\bQ$-rational 
$\Leftrightarrow$ $(2,p)_{2,\bQ(\sqrt{-3})}=0$ 
$\Leftrightarrow$ $p=2,3$, $p\equiv\pm 1\pmod{8}$ $(p\geq 5)$ or 
$p\equiv 2\pmod{3}$ $(p\geq 5)$. 
Note that for $p\geq 5$ 
the condition $p\equiv\pm 1\pmod{8}$ or $p\equiv 2\pmod{3}$ 
is equivalent to $p\not\equiv 13, 19\pmod{24}$.

(iii) For the general case, write $c=\pm p_1^{e_1}\cdots p_r^{e_r}$ where $p_i$ are distinct prime numbers, 
then $L_0$ is $\bQ$-rational 
$\Leftrightarrow$ $(2,c)_{2,\bQ(\sqrt{-3})}=0$ 
$\Leftrightarrow$ $e_i$ is an even integer whenever 
$p_i\equiv 13,19 \pmod{24}$ for any $1\leq i\leq r$.\\

(2) The case where $X_0$ is singular, i.e. $d^2-4bc^2=0$. 

We will construct an example such that 
$X_0$ has only the one (singular) $k$-point and 
$X_0$ is not $k$-rational (thus $X_0$ is not $k$-unirational).

We take $k=\bQ$ and 
$(a,b,c,d,\alpha,\beta)= (3,1,14,28,-1,0)$ with $\alpha^2-a\beta^2=b$. 
Note that, because $d^2-4bc^2=0$, 
it is necessary that $b$ is a square in $k\setminus \{0\}$. 

Substitute these values into Equation (\ref{eq11}).
Then $(\alpha,\beta,\gamma,\delta)=(-1,0,0,0)$ is a solution of
Equation (\ref{eq11}) for these parameters $a, b, c, d$.
The point $(-1,0,0,0)=(-\tfrac{d}{2c},0,0,0)$ becomes 
the unique singular $\bQ$-point on $X_0$.

On the other hand, we have $2c\alpha+d=-28+28=0$. 
Hence it remains to check $(a,2d)_{2,\bQ}=(3,56)_{2,\bQ}\neq 0$ 
by Theorem \ref{t1.8} (viii). 

In the next paragraph, we will show that
$(3,56)_{2,\bQ}=(3,14)_{2, \bQ}\neq 0$; 
thus $L_0$ and $X_0$ are not $\bQ$-rational. 
By Theorem \ref{t1.8}, $X_0$ can have only one (singular) $\bQ$-point. 
To explain these pathologies, note that it arises from 
the fact that $X_0$ fails to be $\bQ$-unirational.

To show that $(3,14)_{2,\bQ}\neq 0$, 
we will prove a more general case: 
For any integer $m\equiv 2\pmod{3}$, $(3,m)_{2,\bQ}\neq 0$. 

Suppose that $(3,m)_{2,\bQ}=0$. 
Then there are integers $u,v,w$ satisfying that 
$3u^2+mv^2=w^2$ with $u,v,w$ not all zero.
We may also assume that $\gcd \{u,v,w \}=1$.
Taking modulo $3$, we get $mv^2\equiv w^2\pmod{3}$.

If $3 \mid w$, then $3 \mid v$.
It follows that $9 \mid 3u^2$. Hence $3 \mid u$.
A contradiction.

If $3 \nmid w$, then $3 \nmid v$.
Thus $m$ is a quadratic residue mod $3$.
Hence the Legendre symbol $(\frac{m}{3})=1$.
But $(\frac{m}{3})=(\frac{2}{3})=-1$ 
because $m\equiv 2\pmod{3}$. 
Another contradiction. Done.
\end{example}
%

\section{Appendix: a geometric interpretation}\label{se4}

After the first version of this paper posted in the arXiv, 
J.-L. Colliot-Th\'el\`ene gave the authors how 
the presence of a smooth $k$-point implies the $k$-rationality 
of $X$ 
when a base field $k$ is perfect infinite with {\rm char} $k\neq 2$. 

In this section, we assume that 
a base field $k$ is perfect infinite with {\rm char} $k\neq 2$. 
The intersection of two quadrics in $\bP^n$ where $n$ is an even integer 
and $n\ge 4$ was studied by Colliot-Th\'el\`ene, Sansuc and Swinnerton-Dyer 
in \cite{CTSSD} by very powerful methods. 
In particular, the theorem \cite[Theorem 2.4]{CTSSD} 
(see Theorem \ref{t3.1} below) is a convenient tool 
to proving the $k$-rationality when a smooth $k$-point exists. 
We will give an alternative proof of the equivalence of 
{\rm (i)}, {\rm (ii)}, {\rm (iii)}, {\rm (iv)}, {\rm (v)} 
of Theorem \ref{t1.8} by using the theorem \cite[Theorem 2.4]{CTSSD}. 

We first explain a terminology used in \cite{CTSSD}. 
Let $\overline{k}$ be the algebraic closure of $k$. 
If $\Pi_1, \Pi_2$ are linear subspaces of $\bP^n_{\overline{k}}$, 
define ${\rm Span}\{\Pi_1, \Pi_2\}$ to be the linear subspace 
spanned projectively by $\Pi_1$ and $\Pi_2$ in $\bP^n_{\overline{k}}$.
The special case of \cite[Theorem 2.4]{CTSSD} when
$n=4$ may be formulated as follows.

\begin{theorem}[Colliot-Th\'el\`ene, Sansuc and Swinnerton-Dyer {\cite[Theorem 2.4]{CTSSD}}] \label{t3.1}
Let $k$ be a perfect infinite field with {\rm char} $k\neq 2$,
$X$ be a projective surface in $\bP^4_{k}$ and
$\overline{X}:= X \times_k \overline{k}$ such that $\overline{X}$
is an integral scheme over $\overline{k}$, pure of dimension two,
and $X$ is an intersection of two quadrics defined over $k$.
Assume that $X$ contains a smooth $k$-point $P_0$ and $\overline{X}$
contains a pair of lines $\Pi_1$ and $\Pi_2$ satisfying that
{\rm (i)} $\Pi_1 \cap \Pi_2 = \emptyset$, $\Pi_1, \Pi_2$ are permuted by
${\rm Gal}(\overline{k}/k)$; and
{\rm (ii)} $P_0 \notin {\rm Span}\{\Pi_1, \Pi_2 \}$,
${\rm Span}\{P_0, \Pi_1 \} \not\subset \overline{X}$ and
${\rm Span}\{P_0, \Pi_2 \} \not\subset \overline{X}$.
Then $X$ is $k$-rational.
\end{theorem}

{\it Alternative proof of 
${\rm (i)}\Leftrightarrow {\rm (ii)}\Leftrightarrow {\rm (iii)}
\Leftrightarrow {\rm (iv)}\Leftrightarrow {\rm (v)}$ 
of Theorem \ref{t1.8} when $k$ is perfect infinite.}

The implications ${\rm (i)}\Rightarrow {\rm (ii)}\Rightarrow {\rm (iii)}$ 
are standard. 
The equivalence ${\rm (iv)}\Leftrightarrow {\rm (v)}$ 
follows from the same argument given in Section \ref{se2}. 
We will give an alternative proof of 
${\rm (iii)}\Rightarrow {\rm (v)}$ and 
${\rm (v)}\Rightarrow {\rm (i)}$. 

\bigskip
{\it Proof of ${\rm (iii)}\Rightarrow {\rm (v)}$ of Theorem \ref{t1.8}.} ---------

Since $X$ is $k$-unirational, there is a dominant rational map 
$f: \bP^m_{k}\dashrightarrow X$ over $k$ 
($m$ is some positive integer). 
Because $k$ is an infinite field, 
the image of $f$ contains infinitely many $k$-points in $X$.

Alternatively, these implications can be found in the proof 
of ${\rm (iii)}\Rightarrow {\rm (v)}$ of Theorem \ref{t1.8}, 
because we have infinitely many choices of $a_0, b_0 \in k$ 
with $h(a_0, b_0) \neq 0$.

\bigskip
{\it Proof of ${\rm (v)}\Rightarrow {\rm (i)}$ of Theorem \ref{t1.8}.} ---------

Let $P_0$ be a smooth $k$-point on $X$. 
Recall that $P_0$ must lie on $X_0$ and hence its $T_5$-coordinate is not zero 
(see the first paragraph of Proof of Theorem \ref{t1.8} in Section \ref{se2}). 
Define two lines $\Pi_1$ and $\Pi_2$ on $\overline{X}$ as follows: 
$\Pi_1$ is defined by the equations $T_1- \sqrt{a}T_2 =T_3- \sqrt{a}T_4=T_5=0$,
$\Pi_2$ is defined by the equations $T_1+\sqrt{a}T_2 =T_3+\sqrt{a}T_4=T_5=0$. 
Clearly $\Pi_1 \cap \Pi_2 = \emptyset$ and $\Pi_1, \Pi_2$ are permuted by 
${\rm Gal}(\overline{k}/k)$. 
Moreover, $P_0 \notin {\rm Span}\{\Pi_1, \Pi_2 \}$, 
because both $\Pi_1$ and $\Pi_2$ lie on the hyperplane $T_5=0$, 
while $P_0$ is a $k$-point on $X$ whose $T_5$-coordinate is not zero.

Finally we will show that ${\rm Span}\{P_0, \Pi_1 \} \not\subset \overline{X}$.
Suppose that ${\rm Span}\{P_0, \Pi_1 \} \subset \overline{X}$. 
We will find a contradiction.

Note that the projective subvarieties $\Pi_1$, ${\rm Span}\{P_0, \Pi_1 \}$, 
$\overline{X}$ in $\bP^4_{\overline{k}}$ are irreducible (integral).

Since $\Pi_1$ is a proper subset of ${\rm Span}\{P_0, \Pi_1 \}$, 
it follows that dim ${\rm Span}\{P_0, \Pi_1 \} >$ dim $\Pi_1 =1$. 
On the other hand, as we suppose that 
${\rm Span}\{P_0, \Pi_1 \} \subset \overline{X}$, 
we have dim ${\rm Span}\{P_0, \Pi_1 \} \le$ dim $\overline{X}=2$. 
It follows that both dimensions of ${\rm Span}\{P_0, \Pi_1 \}$
and $\overline{X}$ are equal to $2$. 
Thus ${\rm Span}\{P_0, \Pi_1 \} = \overline{X}$ 
because $\overline{X}$ is integral. 
Consequently $\overline{X}$ is also a linear variety 
(a linear space in $\bP^4_{\overline{k}}$).

The degree of $\overline{X}$ is $1$ because it is a linear variety \cite[page 52]{Ha}. However, $\overline{X}$ is of degree $4$, because it is the intersection of two quadrics. This leads to a contradiction.

By the same way, ${\rm Span}\{P_0, \Pi_2 \} \not\subset \overline{X}$.

Now we may apply Theorem \ref{t3.1} to conclude that $X$ is $k$-rational.
\qed


\end{document}